\documentclass[12pt]{conm-p-l}
\usepackage{amssymb, amstext, amscd, amsmath}
\makeatletter
\def\@cite#1#2{{\m@th\upshape\bfseries%
[{#1\if@tempswa{\m@th\upshape\mdseries, #2}\fi}]}}
\makeatother
\theoremstyle{plain}
\newtheorem{thm}{Theorem}[section]
\newtheorem{cor}[thm]{Corollary}

\theoremstyle{definition}

\newtheorem{eg}[thm]{Example}

\newcommand{\bC}{{\mathbb{C}}}

\newcommand{\bD}{{\mathbb{D}}}
\newcommand{\bK}{{\mathbb{K}}}
\newcommand{\bN}{{\mathbb{N}}}

\newcommand{\bR}{{\mathbb{R}}}
\newcommand{\bT}{{\mathbb{T}}}
\newcommand{\bZ}{{\mathbb{Z}}}

  \newcommand{\A}{{\mathcal{A}}}
  \newcommand{\B}{{\mathcal{B}}}

\renewcommand{\H}{{\mathcal{H}}}

  \newcommand{\K}{{\mathcal{K}}}

\renewcommand{\O}{{\mathcal{O}}}

  \newcommand{\U}{{\mathcal{U}}}

\renewcommand{\phi}{\varphi}

\newcommand{\rC}{{\mathrm{C}}}

\newcommand{\fA}{{\mathfrak{A}}}
\newcommand{\fB}{{\mathfrak{B}}}

\newcommand{\fM}{{\mathfrak{M}}}

\newcommand{\fs}{{\mathfrak{s}}}

\newcommand{\bi}{{\mathbf{i}}}
\newcommand{\bn}{{\mathbf{n}}}
\newcommand{\bx}{{\mathbf{x}}}

\newcommand{\AND}{\text{ and }}
\newcommand{\FOR}{\text{ for }}

\newcommand{\qand}{\quad\text{and}\quad}

\newcommand{\qfor}{\quad\text{for}\quad}
\newcommand{\qforal}{\quad\text{for all}\quad}

\newcommand{\ad}{\operatorname{ad}}
\newcommand{\Alg}{\operatorname{Alg}}

\newcommand{\id}{{\operatorname{id}}}

\newcommand{\ca}{\mathrm{C}^*}

\newcommand{\cenv}{\mathrm{C}^*_{\text{env}}}
\newcommand{\ol}{\overline}

\newcommand{\Fn}{\mathbb{F}_n^+}

\newcommand{\Fock}{\ell^2(\Fn)}

\newcommand{\mt}{\varnothing}
\newcommand{\ltwo}{\ell^2}

\newcommand{\diag}{\operatorname{diag}}

\begin{document}

%%%%%%%%%%%%%%%%%%%%%%%%%%%%%%%%
\title[Nonself-adjoint operator algebras]%
{Nonself-adjoint operator algebras\\ for dynamical systems}
\thanks{}

\author[K.R.Davidson]{Kenneth R. Davidson}
\address{Pure Math.\ Dept.\\U. Waterloo\\Waterloo, ON\;
N2L--3G1\\CANADA}
\email{krdavids@uwaterloo.ca}

\author[E.G. Katsoulis]{Elias~G.~Katsoulis}
\address{Dept.\ Mathematics\\East Carolina U.\\
Greenville, NC 27858\\USA}
\email{KatsoulisE@mail.ecu.edu}
\begin{abstract}
This paper is a survey of our recent work on operator algebras
associated to dynamical systems that lead to classification results
for the systems in terms of algebraic invariants of the operator algebras.
\end{abstract}

\thanks{First author partially supported by an NSERC grant.}
\thanks{Second author partially supported by
a summer grant from ECU}
\date{}
\maketitle

%%%%%%%%%%%%%%%%%%%%%%%%%%%%%%%%
\section{Introduction}
%%%%%%%%%%%%%%%%%%%%%%%%%%%%%%%%

There is a long history of building operator algebras from dynamical systems
going back to von Neumann's construction of a group von Neumann algebra.
The use of nonself-adjoint operator algebras is more recent, but still goes back
40 years to Arveson's paper \cite{Arv}.  His algebra was closed in the weak operator
topology.  We will mostly be interested in norm closed (generally separable) algebras,
but these algebras are the same vintage, beginning with work of
Arveson and Josephson \cite{ArvJ}.

For the most part, we consider a dynamical system to be a locally compact
Hausdorff space $X$ together with one or more proper maps $\sigma_i$, $1 \le i \le n$,
of $X$ into itself, which we write as $(X,\sigma)$.
A proper map $\sigma$ of $X$ into itself induces an endomorphism
$\alpha(f) = f \circ\sigma$ of $\rC_0(X)$.  So more generally, we could consider
an arbitrary semigroup of endomorphisms of an arbitrary algebra.
This is too general a setting, but we will consider some variants of this kind.
In particular, we will consider a C* dynamical system $(\fA,\alpha)$, where $\fA$
is a C*-algebra with a single automorphism $\alpha$.

Two systems $(X,\sigma)$ and $(Y,\tau)$ are said to be \textit{conjugate} provided that
there is a homeomorphism $\gamma$ of $X$ onto $Y$ so that
$\tau_i\gamma = \gamma\sigma_i$ for $1 \le i \le n$.
Similar notions may defined for endomorphisms of other operator algebras as well.
The main question is whether the operator algebra encodes enough information
intrinsically so that the dynamical system can be recovered up to conjugacy
or some similar relation between dynamical systems.

Our operator algebras can be considered in the abstract, in the spirit of the
Blecher--Ruan--Sinclair Theorem \cite{BRS}.
Generally they are defined by a universal property dealing with representations.
So the question of an sufficient family of such representations is moot.
See \cite{Pau, BL} for more on this viewpoint.

%%%%%%%%%%%%%%%%%%%%%%%%%%%%%%%%
\section{Semicrossed Products}
%%%%%%%%%%%%%%%%%%%%%%%%%%%%%%%%

Arveson \cite{Arv} and Arveson--Josephson \cite{ArvJ} were the first to
consider nonself-adjoint operator algebras constructed from dynamical systems.
They considered a homeomorphism $\sigma$ of a compact Hausdorff space $X$
which preserves a probability measure $\mu$ satisfying $\mu(O)>0$ for every non-empty
open set $O\subset X$.  Given such a system, one can define a representation of
$\rC(X)$ on $L^2(\mu)$ as multiplication operators, and define a unitary
$Uf = f \circ\sigma$.  Let $\A(X,\sigma)$ denote the norm-closed subalgebra
of $\B(L^2(\mu))$ generated by $U$ and the multiplication operators.
Their main result is that, under a further technical condition, that the
algebraic structure of the algebra determines the dynamical system:

\begin{thm}[Arveson--Josephson \cite{ArvJ}]
Suppose that for $i=1,2$,
$\sigma_i$ is a homeomorphism of a compact Hausdorff space $X_i$
which preserves a probability measure $\mu_i$ satisfying $\mu_i(O)>0$ for every non-empty
open set $O$ of $X_i$.   Moreover suppose that the set of periodic points has measure $0$.
Then the following are equivalent:
\begin{enumerate}
\item $(X_1,\sigma_1)$ and $(X_2,\sigma_2)$ are conjugate.
\item $\A(X_1,\sigma_1)$ and $\A(X_2,\sigma_2)$ are isometrically  isomorphic.
\item  $\A(X_1,\sigma_1)$ and $\A(X_2,\sigma_2)$ are algebraically isomorphic.
\end{enumerate}
\end{thm}

In 1985, Peters \cite{Pet} introduced an abstract construction of an operator algebra
associated to the system $(X,\sigma)$.  He did not require the
map $\sigma$ to be a homeomorphism, only that it be a proper map so that it induces an
endomorphism $\alpha(f) = f \circ\sigma$ of $\rC_0(X)$.
Moreover, he does not require an invariant measure.

Suppose that $X$ is compact.
The idea is to construct a universal operator algebra which contains $\rC(X)$
as a C*-subalgebra, that is generated by $\rC(X)$ and a single isometry $\fs$
which encodes the dynamics via the \textit{covariance relation}
\[ f \fs = \fs (f\circ\sigma) \qforal f \in \rC_0(X) .\]
Indeed, consider a covariant representation of this system as a pair $(\rho,S)$
where $\rho $ is a $*$-representation of $\rC(X)$ and $S$ is an isometry satisfying
$\rho(f)S = S \rho(f\circ\sigma)$ for all $f \in \rC(X)$.
It is not difficult to show that there is a unique operator algebra
$\rC(X) \times_\alpha\bZ^+$ satisfying these properties together with the
additional property that for every covariant representation $(\rho,S)$,
there is a completely contractive representation $\pi$ of $\fA \times_\alpha\bZ^+$
such that $\pi|_{\rC(X)} = \rho $ and $\pi(\fs) = S$.
This is called the \textit{semi-crossed product algebra}.

In the non-compact case, we do not include $\fs$ in the algebra.
$\rC_0(X) \times_\alpha\bZ^+$ is the universal algebra generated by
$\rC_0(X)$ and the elements $\fs f$ for $f \in \rC_0(X)$.  The isometry $\fs$
still belongs to the multiplier algebra however, and the covariance relations
make sense if we multiply on the right by a function $g\in\rC_0(X)$.
It has the property that for every covariant representation $(\rho,S)$,
there is a completely contractive representation $\pi$ of $\fA \times_\alpha\bZ^+$
such that $\pi|_{\rC(X)} = \rho $ and $\pi(\fs f) = S\rho(f)$ for every $f\in\rC_0(X)$.

Every dynamical system $(X,\sigma)$ has a natural family of covariant representations
called \textit{orbit representations}.  For any point $x\in X$, build a
representation $\pi_x$ on $\ltwo$ by
\[ \pi_x(f) = \diag( f(x), f(\sigma(x)), f(\sigma ^2(x)), \dots) \]
and $\pi_x(\fs) = S_+$, the unilateral shift.
Peters shows that $\Pi = \sum_{x\in X} \pi_x$ is a faithful completely
isometric representation of the semicrossed product $\rC_0(X) \times_\alpha\bZ^+$.
So it has a rather concrete spatial representation.

Peters was able to significantly extend the Arveson--Josephson Theorem.  He still assumes
that $X$ is compact, but works with an arbitrary continuous map $\sigma$ of $X$ into itself,
subject to the assumption that there are no fixed points.

\begin{thm}[Peters \cite{Pet}]
Let $(X_i,\sigma_i)$ be compact dynamical systems with no fixed points.
Then the following are equivalent:
\begin{enumerate}
\item $(X_1,\sigma_1)$ and $(X_2,\sigma_2)$ are conjugate.
\item $\rC(X_1) \times_{\sigma_1}\bZ^+$ and $\rC(X_2) \times_{\sigma_2}\bZ^+$
 are completely isometrically  isomorphic.
\item  $\rC(X_1) \times_{\sigma_1}\bZ^+$ and $\rC(X_2) \times_{\sigma_2}\bZ^+$
 are algebraically isomorphic.
\end{enumerate}
\end{thm}

In 1988, Hadwin and Hoover \cite{HadH} considered a much more general class of
\textit{conjugacy algebras} which contain $\rC(X)$ and an element $\fs$
satisfying the covariance relations.
They do not even require the algebras to be closed.
Their methods weakened the condition on fixed points
to saying: \textit{$\{x \in X_1: \sigma_1^2(x) = \sigma_1(x) \ne x \}$ has no interior.}  With this
hypothesis, they reach the same conclusion.  Moreover, they show that the algebraic
isomorphism of any single conjugacy algebra for the system $(X_1,\sigma_1)$ with any
conjugacy algebra for the system $(X_2,\sigma_2)$ is enough to imply conjugacy
of the systems, and hence completely isometric isomorphism of the semicrossed products.

Another result along these lines is due to Power \cite{Pow}.  He loosens the condition
on $X$ to be locally compact, but requires the map $\sigma$ to be a homeomorphism.
The key advance is that he recovers the Arveson--Josephson conclusions without any
hypothesis on fixed points.

In 2006, the authors removed all extraneous hypotheses on the dynamical system
to obtain:

\begin{thm}[Davidson--Katsoulis \cite{DKisotop}]
Let $X_i$ be locally compact Hausdorff spaces, and let $\sigma_i$ be a proper continuous
map of $X_i$ into itself.  Then the following are equivalent:
\begin{enumerate}
\item $(X_1,\sigma_1)$ and $(X_2,\sigma_2)$ are conjugate.
\item $\rC(X_1) \times_{\sigma_1}\bZ^+$ and $\rC(X_2) \times_{\sigma_2}\bZ^+$
 are completely isometrically  isomorphic.
\item  $\rC(X_1) \times_{\sigma_1}\bZ^+$ and $\rC(X_2) \times_{\sigma_2}\bZ^+$
 are algebraically isomorphic.
\end{enumerate}
\end{thm}

We will discuss some of the ideas that go into the proof.
Hadwin and Hoover study the characters of $\rC(X) \times_{\sigma}\bZ^+$.
The restriction of a character $\theta$ to $\rC_0(X)$ will be a point evaluation
map $\delta_x$ for some $x \in X$.
The rest of $\theta$ is determined by $\theta(\fs) = z$.
As $\theta$ is contractive, one has $|z|\le1$.
Call this character $\theta_{x,z}$.
It is not difficult to show that $\theta_{x,0}$ always exists,
and that this is the only possibility if $\sigma(x) \ne x$.
When $\sigma(x)=x$, then $\theta_{x,z}$ exists for all $z \in \ol{\bD}$.
We write $\fM$ for the character space, and $\fM_x$ for those characters
which restrict to $\delta_x$ on $\rC_0(X)$.

We introduce the notion of an \textit{analytic set} in the character space.
This is the image of a continuous non-constant map $F$ of a domain
$\Omega \subset \bC$ into $\fM$ such that $f_a(w) = F(w)(a)$ is analytic on $\Omega$
for all $a\in \rC(X) \times_{\sigma}\bZ^+$.
A subset of $\fM$ is a \textit{maximal analytic set} if it is an analytic set
which is not properly contained in any other.
A crucial observation is that $\{ \theta_{x,z} : z \in \bD \}$ are maximal analytic sets
for every fixed point $x = \sigma(x)$, and there are no others.

The other important device is the notion of a nest representation.
A representation of an algebra into $\B(\H)$ is said to be a \textit{nest representation}
provided that the lattice of invariant subspaces of its range is a nest (complete chain) of subspaces.
Nest representations and their kernels were originally introduced
by Lamoureux \cite{Lam} as a generalization for primitive ideals.
The second author, in joint work with Peters \cite{KaPet2} and Kribs \cite{KK}, applied the
concept of a nest representation to the classification theory for
non-selfadjoint operator algebras, including
graph algebras \cite{KK} and limit algebras \cite{KaPet2}. (See also \cite{Sol}.)
For our purposes, nest
representations into the $2\times 2$ upper triangular
matrices are sufficient.

In our case, such a $2\times2$ nest representation of $\rC(X) \times_{\sigma}\bZ^+$
will have diagonal entries which are characters $\theta_{x,z}$ and $\theta_{y,w}$,
and a $1,2$ entry which is linearly independent of the diagonal (in order that there
be no other invariant subspace).  The dynamical system is recovered by showing
that this implies that $\sigma(x) = y$.  A complication occurs when $\sigma(y)=y$.
This is the case that forced Hadwin and Hoover to include an additional hypothesis.
We circumvent that by considering an analytic family of nest representations.

%%%%%%%%%%%%%%%%%%%%%%%%%%%%%%%%
\section{Function Algebra Systems}
%%%%%%%%%%%%%%%%%%%%%%%%%%%%%%%%

Peter's construction was actually developed in greater generality.
If $\fA$ is any operator algebra and $\alpha$ is a completely contractive
endomorphism, one defines the semicrossed product $\fA \times_\alpha \bZ^+$
in an analogous manner.  This is the universal operator algebra which
is generated by a completely isometric copy of $\fA$ and an isometry $\fs$
satisfying the covariance relation
\[ a \fs = \fs \alpha(a) \qforal a \in \fA .\]
Given any covariant representation $(\rho,S)$ consisting of a
completely contractive representation $\rho$ of $\fA$ and an isometry $S$,
there is a completely contractive representation $\pi$ of $\fA \times_\alpha \bZ^+$
such that $\pi|_\fA = \rho$ and $\pi(\fs) = S$.
Note that when $\fA$ is a C*-algebra, the completely contractive representations are
precisely the $*$-representations.

There is a natural analogue of the orbit representations for this semicrossed product.
Starting with a completely contractive representation $\rho$ of $\fA$ on $\H$,
form a representation on $\ltwo(\H)$ by setting
\[ \pi_\rho(f) = \diag( \rho(a), \rho(\alpha(a)), \rho(\alpha^2(a)),\dots)
  \AND \pi_\rho(\fs) = S_+\otimes I_\H .
\]
Again the direct sum of these representations yields a faithful completely isometric representation of $\fA \times_\alpha \bZ^+$.

To illustrate how the ideas of the previous section can be applied
in another setting, we consider certain function algebras.
Let $G$ be a Cauchy domain in $\bC$, namely a bounded open subset
such that $\partial G$ consists of a finite union of disjoint Jordan curves.
Let $K = \ol{G}$.  Then $A(K)$ is the algebra of continuous functions on $K$
which are analytic on $G$.
Suppose that $\sigma \in A(K)$ is a non-constant function such that
$\sigma(K) \subset K$.  Then $\alpha(f) = f\circ\sigma$ yields a completely
contractive endomorphism of $A(K)$.
Thus we can form the semicrossed product $A(K) \times_\sigma \bZ_+$.

We say that two such maps $\sigma_i:K_i \to K_i$ are \textit{analytically conjugate}
if there is a biholomorphic map of $G_1$ onto $G_2$ such that
$\sigma_2\gamma = \gamma\sigma_1$.
(Any biholomorphic map between bounded domains of $\bC$ extends
to a homeomorphism of the closures.)

\begin{eg}[Buske--Peters \cite{BP}]
Let $\sigma$ be an elliptic Mobius map of $\ol{\bD}$ onto itself
(i.e. $\sigma$ has a single fixed point, and it lies in the open disk $\bD$).
Then $\sigma$ is analytically conjugate to to a rotation $\eta(z) = e^{2\pi \theta i} z$.
They show that $A(\bD) \times_\sigma \bZ_+$ is isometrically isomorphic to
the subalgebra $\Alg\{U,V\}$ of the rotation C*-algebra $\A_\theta$.
It follows that there is an automorphism interchanging $U$ and $V$ which
has the effect of flipping the automorphism $\sigma$ to $\sigma^{-1}$.
So  $A(\bD) \times_\sigma \bZ_+$ and $A(\bD) \times_{\sigma^{-1}} \bZ_+$
are completely isometrically isomorphic.
In particular, the dynamics cannot be completely recovered from
the semicrossed product in this case.
\end{eg}

It turns out that this is the only thing which can complicate matters.

\begin{thm}[Davidson--Katsoulis \cite{DKisotop}]
For $i=1,2$, let $G_i$ be Cauchy domains in $\bC$ with $K_i = \ol{G_i}$.
Then the following are equivalent:
\begin{enumerate}
\item $A(\bK_1) \times_{\sigma_1} \bZ_+$ and $A(\bK_2) \times_{\sigma_2} \bZ_+$
are   algebraically isomorphic.
\item $A(\bK_1) \times_{\sigma_1} \bZ_+$ and $A(\bK_2) \times_{\sigma_2} \bZ_+$
are completely isometrically isomorphic.
\item {\em (i)} $\sigma_1$ and $\sigma_2$ are analytically conjugate, or\\
{\em (ii)} $\bK_i$ are simply connected, $\sigma_i$ are homeomorphisms with a
unique fixed point which lies in the interior $G_i$, and
$\sigma_2$ is analytically conjugate to $\sigma_1^{-1}$.
\end{enumerate}
\end{thm}

When one applies the same analysis to $A(K)$, one finds an additional maximal
analytic set in the character space, namely $\{\theta_{x,0} : x \in G \}$.
This causes no difficulties when $G$ is not simply connected because
it must be identified with the corresponding non-simply connected analytic
set of the other algebra.  However when $G$ is conformally equivalent to the disk,
there can be some interplay---and this is precisely what occurs in the Buske--Peters example.
Only the elliptic case is complicated, as otherwise the fixed point(s) are on the
boundary, and topological considerations show that the new analytic disk
is uniquely determined.  In the elliptic case, one needs to work harder to
show that an algebraic isomorphism is automatically continuous.
This makes use of ideas from a paper of Sinclair \cite{Sinclair}.

Finally, we mention that the case of $\sigma = \id$ actually requires
a special result from function theory.  In this case,
$A(K) \times_\sigma \bZ_+ \simeq A(K \times \ol{\bD})$.
We use the fact that any biholomorphic map of one product
space onto another actually decomposes as a product map \cite{Lig}.

%%%%%%%%%%%%%%%%%%%%%%%%%%%%%%%%
\section{C*-dynamical Systems}
%%%%%%%%%%%%%%%%%%%%%%%%%%%%%%%%

In this section, we will consider the (nonself-adjoint) semicrossed product
constructed from a pair $(\fA,\alpha)$ where $\fA$ is a C*-algebra and
$\alpha$ is a $*$-automorphism.
In this case, $\fA \times_\alpha \bZ_+$ is the subalgebra of the
crossed product C*-algebera $\fA \times_\alpha \bZ$ generated by $\fA$ and
the unitary $U$ implementing $\alpha$.

Two C*-dynamical systems $(\fA,\alpha)$ and $(\fB,\beta)$
are \textit{conjugate} if there is a $*$-isomorphism $\gamma$ of $\fA$ onto $\fB$
so that $\beta\gamma = \gamma \alpha $.
In the non-commutative setting, this is too strong because one can modify
$\alpha$ by an inner automorphism of $\fA$ without affecting the algebra structure
of the crossed product.
Therefore we say that these systems are \textit{outer conjugate}
if there is a unitary element $u\in\fA$ so that $\beta\gamma = \gamma \ad_u\alpha$.

Automorphisms of C*-algebras have been widely studied.
One useful tool is the \textit{Connes spectrum}.  Let $\H^\alpha(\fA)$
denote the set of all hereditary sub-C*-algebras of $\fA$ which are
$\alpha$-invariant.  Define a subset of the circle $\bT$ by
\[ \Gamma(\alpha) = \cap _{\fB \in \H^\alpha(\fA)} \sigma(\alpha|_\fB) .\]
This is in some sense analogous to the Weyl spectrum of a bounded operator.
The following result of Olesen and Pedersen \cite{OPed} characterizes the
important special case of full Connes spectrum, $\Gamma(\alpha) = \bT$.

We write $\hat\fA$ denote the spectrum of $\fA$, and let $\hat\alpha$ be
the induced action of $\alpha$ on $\hat\fA$ by $\hat\alpha([\pi]) = [\pi\alpha]$.

\begin{thm}[Olesen--Pedersen \cite{OPed}]
Let $\fA$ be a separable C*-algebra, and fix an automorphism $\alpha$ of $\fA$.
Then the following are equivalent:
\begin{enumerate}
\item $\alpha$ has full Connes spectrum, $\Gamma(\alpha) = \bT$.
\item there is a dense subset $\Delta$ of the spectrum $\hat{\fA}$ which
is $\hat\alpha$-inv\-ar\-iant on which $\hat\alpha$ acts freely.
\item $\{ [\pi] \in \hat\fA : \hat\alpha^n([\pi]) = [\pi] \}$ has no interior for all $n\ge1$.
\end{enumerate}
\end{thm}

This was used by Muhly and Solel  to show the following:

\begin{thm}[Muhly--Solel \cite{MS}]
Suppose that $(\fA,\alpha)$ and $(\fB,\beta)$ are C*-dynamical systems
such that $\Gamma(\alpha) = \bT$.
Then $\fA \times_\alpha \bZ_+$ and $ \fB \times_\beta \bZ_+$ are completely
isometrically isomorphic
if and only if $(\fA,\alpha)$ and $(\fB,\beta)$ are outer conjugate.
\end{thm}

We say that an automorphism $\alpha$ is \textit{universally weakly inner}
with respect to irreducible representations if for every irreducible representation $\pi$
of $\fA$, there  exists a unitary $W \in \pi(A)''$ so that $\pi\alpha(A) = W\pi(A) W^*$.
Kishimoto \cite{Kis} proves the remarkable result that if $\fA$ is a simple
separable C*-algebra, then
every universally weakly inner automorphism is actually inner.

We were able to avoid the condition on Connes spectrum by modifying the
arguments used in the commutative case.

\begin{thm}[Davidson--Katsoulis \cite{DKsemicross}]
Let $(\fA, \alpha)$ and $(\fB, \beta)$ be C*-dynamical systems, and
assume that the semicrossed products $\fA \times_{\alpha} \bZ^{+}$ and
$\fB \times_{\beta} \bZ^{+}$ are completely isometrically isomorphic. Then
there exists a C*-isomorphism $\gamma: \fA\to\fB$
so that $\alpha\circ\gamma^{-1}\circ \beta^{-1}\circ \gamma$ is
universally weakly inner with respect to irreducible representations.
\end{thm}

Then using Kishimoto's result, we obtain that

\begin{cor}
If $\fA$ is a separable simple C*-algebra, then
$\fA \times_\alpha \bZ_+$ and $\fB \times_\beta \bZ_+$ are isometrically isomorphic
if and only if $(\fA,\alpha)$ and $(\fB,\beta)$ are outer conjugate.
\end{cor}

Our methods yield a new proof of the Muhly--Solel theorem.  They suggest that
isometric isomorphism of the semicrossed product may imply outer
conjugacy under quite general hypotheses.

%%%%%%%%%%%%%%%%%%%%%%%%%%%%%%%%
\section{Multivariable dynamical systems}
%%%%%%%%%%%%%%%%%%%%%%%%%%%%%%%%

A multivariable dynamical system is a locally compact Hausdorff space $X$
together with a collection $\sigma_1,\dots,\sigma_n$ of proper continuous maps
of $X$ into itself.  We will seek an appropriate analogue of the semicrossed product.
To this end, we seek an operator algebra which contains $\rC_0(X)$
and operators $\fs_1,\dots,\fs_n$ satisfying the covariance relations
\[  f \fs_i = \fs_i (f \circ \sigma_i) \qfor 1 \le i \le n \AND f \in \rC_0(X) . \]
Again, in the non-compact case, we do not include the $\fs_i$ in the algebra,
but do include the elements $\fs_i g$.

Since we do not require any relations between the maps, the natural
semigroup that arises is the free semigroup $\Fn$ of all words in an
alphabet of $n$ letters (including the empty word).  If $w=i_1\dots i_k$
is an element of $\Fn$, we write $\fs_w = \fs_{i_1}\dots \fs_{i_k}$.
Similarly, we write $\sigma_w = \sigma_{i_1} \circ \dots \circ \sigma_{i_k}$.

The issue of norming the elements $\fs_i$ leads to a couple of natural choices
that yield different algebras.
The simplest condition is just to insist that each $\|\fs_i\| \le 1$.
The universal algebra subject to this constraint will be called the
semi-crossed product, denoted $\rC_0(X) \times_\sigma \Fn$.
This has the universal property that given $(\rho,S_1,\dots,S_n)$, where $\rho$ is a
$*$-representation of $\rC_0(X)$ on a Hilbert space $\H$, and $S_1,\dots,S_n$
are $n$ contractions on $\H$ satisfying the covariance relations,
then there is a completely contractive representation $\pi$ of
$\rC_0(X) \times_\sigma \Fn$ such that $\pi|_{\rC_0(X)}=\rho$
and $\pi(\fs_i f) = S_i \rho(f)$ for $1 \le i \le n$ and $f\in\rC_0(X)$.

The other reasonable option is to require that $S = \big[ S_1\ \dots\ S_n \big]$
be a row contraction (as an operator from $\H^{(n)}$ to $\H$).
We call the universal algebra obtained in this manner the
\textit{tensor algebra} $\A(X,\sigma)$.

You may notice that we did not require the operators to be isometries.
However such a requirement would make no difference.
The reason is that there is a dilation theorem showing that any
contractive covariant representation of $(X,\sigma)$ dilates to one
in which each $S_i$ is an isometry; and each row contractive
covariant representation dilates to a row isometric covariant
representation.  One obvious advantage of this formulation is
that we now know more about the (completely contractive) representations
of these algebras.

The tensor algebra turns out to be more tractable in general.
One reason is that there is a natural analogue of the orbit
representations.  Fix $x\in X$ and build a representation on
Fock space $\Fock$, with orthonormal basis $\{\xi_w : w \in \Fn \}$, by setting
\begin{align*}
\pi_x(f) \xi_w &= f(\sigma_w(x)) \xi_w\\
\pi_x(\fs_i) \xi_w &= \xi_{iw} \qquad\qquad\quad \FOR w\in\Fn.
\end{align*}
As in the one variable case, the direct sum of all orbit representations
yields a faithful, completely isometric representation of $\A(X,\sigma)$.
The semicrossed product does not appear to have a nice family of
norming representations that can be explicitly described.

In either case, the universal property leads to the existence of
gauge automorphisms, namely automorphisms $\gamma_z$ for $z\in\bT$
of our algebra with $\gamma_z|_{\rC_0(X)} = \id$ and $\gamma_z(\fs_i f) = z\fs_i f$.
In the standard manner, integration over the unit circle yields a completely
contractive expectation onto $\rC_0(X)$.
This leads to an automatic continuity result that any isomorphism of
$\A(X, \sigma)$ onto $\A(Y , \tau)$ or of $\rC_0(X) \times_\sigma \Fn$
onto $\rC_0(Y) \times_\tau \Fn$ is automatically norm continuous.

There are no labels on our maps, or on the isometries $\fs_i$.
So an isomorphism can permute the maps arbitrarily.
It is less obvious, but still true, that in some circumstances,
one can change from one permutation to another.
This leads to our definition of \textit{piecewise conjugate} systems.
Say that $(X,\sigma)$ is piecewise conjugate to $(Y,\tau)$
if there is a homeomorphism $\gamma: X \to Y$
and an open cover $\{ \O_\alpha : \alpha \in S_n \}$ of $X$
so that
\[
 \tau_i \gamma|_{\O_\alpha} = \gamma \sigma_{\alpha(i)}|_{\O_\alpha}
 \qfor \alpha \in S_n .
\]

To appreciate this notion, consider two maps $\sigma_1$ and $\sigma_2$
that map $[0,1]$ into itself and coincide on an interval $(a,b)$.
Then one can construct $\tau_1$ which agrees with $\sigma_1$ on $[0,b)$ and
with $\sigma_2$ on $(a,1]$; and similarly $\tau_2$ agrees with $\sigma_2$
on $[0,b)$ and with $\sigma_1$ on $(a,1]$.  Then $([0,1],\sigma_1,\sigma_2)$
and $([0,1],\tau_1,\tau_2)$ are piecewise conjugate.
On the other hand, if the two maps only agree at a point $\{a\}$,
then one can still define the maps $\tau_1$ and $\tau_2$ as above, but the
new system will not be piecewise conjugate because there is no neighbourhood
of $a$ on which we can match up the two pairs of functions.

This appears to be a new notion in dynamics.  There is a parallel with
the \textit{full group} introduced by Dye \cite{Dye} in his analysis of
group actions on von Neumann algebras.  One begins with a group $\{\alpha_g : g \in G\}$
of measure preserving automorphisms of a measure space $\fM$.  Dye considers all
automorphisms $\alpha$ which are pieced together by a countable partition
of the space into measureable sets $P_g$ such that $\alpha_g^{-1}\alpha$ is the
identity on $P_g \fM$.  So the notion of mixing and matching maps occurs here,
but in a rather different context.

The main result of our paper \cite{DKmultdyn} is the following:

\begin{thm}[Davidson--Katsoulis \cite{DKmultdyn}]
Let $(X,\sigma)$ and $(Y,\tau)$ be two multivariable dynamical systems.
If there is an algebra isomorphism of $\A(X,\sigma)$ onto $\A(Y,\tau)$ or
of $\rC_0(X) \times_\sigma \Fn$ onto $\rC_0(Y) \times_\tau \Fn$,
then $(X,\sigma)$ and $(Y,\tau)$  are piecewise conjugate.
\end{thm}

The proof follows the ideas of the $n=1$ case, but non-trivial
complications arise.  In particular, one must be able to count
the number of maps in the system which send a point $x$ to a point $y$.
The key is again an analytic structure on the set of nest representations.
The ability to count the number of maps relies on the well-known, but
non-trivial, fact from several complex variable theory \cite{FrG} that
the zero set of an analytic function mapping $\bC^k$ into $\bC^l$,
for $l<k$, has no isolated points.

In the case $n=1$, the converse direction was trivial.  But here there
are difficult issues about how to intertwine the isometries to form
the new ones.  This appears to be possible in the case of the tensor algebra.
We have no idea how this could be accomplished in the semicrossed
product case.  A partial converse is the following:

\begin{thm}[Davidson--Katsoulis]
Suppose that at least one of the following holds:
\begin{itemize}
\item $n \le 3$, or
\item $X$ has covering dimension at most 1, or
\item $\{ x : |\sigma(x)| < n \}$ has no interior.
\end{itemize}
Then the following are equivalent:
\begin{enumerate}
\item $(X,\sigma)$ and $(Y,\tau)$ are piecewise conjugate.
\item $\A(X,\sigma)$ and $\A(Y,\tau)$ are algebraically isomorphic.
\item $\A(X,\sigma)$ and $\A(Y,\tau)$ are completely isometrically isomorphic.
\end{enumerate}
\end{thm}

We conjecture that the converse holds in complete generality.
It was reduced to a technical conjecture about the existence of a nice
map from the polytope with vertices indexed by $S_n$ into the unitary
group $\U_n$ which takes the vertices to the corresponding permutation
matrix, and satisfies some strict compatibility conditions on the various faces.
Chris Ramsey, a student at the University of Waterloo, has been making
progress on this conjecture.

It would be interesting to study the ideal structure of these algebras.
Peters \cite{Pet2} has made progress on this in the case $n=1$.
A natural test question is to determine when the algebra is semisimple,
and more generally, to identify the radical.  This has been accomplished
in the case $n=1$ by  Donsig,  Katavolos and Manoussos \cite{DKM},
with earlier work by Muhly \cite{MuRad}.
We can answer the question about semisimplicity, but have little to say
about the radical.

Call an open subset $U \subset X$ a \textit{$(u,v)$--wandering set} if
\[
  \sigma_{u w v}^{-1}(U) \cap U = \mt   \qforal w \in \Fn .
\]
A generalized wandering set is a $(u,v)$-wandering set for some pair $(u,v)$.
If there are no wandering sets, then necessarily each $\sigma_i$ is surjective.

Wandering sets have a parallel notion of recurrence.
Say that $x \in X$ is \textit{$(u,v)$--recurrent} if for every open set $U \ni x$,
there is some $w\in\Fn$ so that $\sigma_{uwv}(x) \in U$.
In the metrizable case, there are no non-empty generalized wandering sets
if and only if the $(u,v)$--recurrent points are dense for all pairs $(u,v)$.  This is
in turn equivalent to the surjectivity of each $\sigma_i$ and the density of the
$(\mt,v)$--recurrent points for each $v\in\Fn$.

\begin{thm}[Davidson--Katsoulis]
The following are equivalent:
\begin{enumerate}
\item $\A( X , \sigma)$  is semisimple.
\item $\rC_0(X) \times_\sigma \Fn$  is semisimple.
\item There are no non-empty generalized wandering sets.
\end{enumerate}
\end{thm}

%%%%%%%%%%%%%%%%%%%%%%%%%%%%%%%%
\section{C*-envelopes}
%%%%%%%%%%%%%%%%%%%%%%%%%%%%%%%%

Arveson's seminal paper \cite{Arv1} proposes that to study a nonself-adjoint
operator algebra, there should be a canonical minimal C*-algebra that contains it
(completely isometrically).  This C*-algebra, $\cenv(\A)$, is the proposed analogue of the
Shilov boundary in the function algebra case, and is called
the \textit{C*-envelope} of $\A$.  Let $j_0$ be
the completely isometric imbedding of $\A$ into $\cenv(\A)$.  Then $\cenv(\A)$  is determined by the
universal property that whenever $j$ is a completely isometric isomorphism of $\A$ into
another C*-algebra $\ca(j(\A))$, there exists a $*$-homomorphism $\pi$ of
$\ca(j(\A))$ onto $\cenv(\A)$ such that $\pi j = j_0$.

Unlike the other universal constructions mentioned in this paper, it is not
at all apparent that the C*-envelope exists.  Arveson constructed it for a
large family of examples, but left the existence in general as a conjecture.
This was verified a decade later by Hamana \cite{Ham}.
A new proof was found a few years ago by Dritschel and McCullough \cite{DMc}.
Their proof is based on the notion of a \textit{maximal dilation}.
A representation $\rho$ of $\A$ on $\H$ is maximal if any dilation of $\rho$
to a completely contractive representation $\pi$ on a larger space $\K$
(meaning that $\rho(A) = P_\H \pi(A)|_\H$) has the form $\pi = \rho \oplus \pi'$
on $\K = \H \oplus \H^\perp$.  It is not particularly difficult to show that
any completely contractive representation can be dilated to a maximal one.
The point is that maximal representations extend to $*$-representations
of the enveloping C*-algebra of $\A$ and factor through the C*-envelope.
Some of these ideas were already known due to work of Muhly and Solel \cite{MSbound}.
This new proof provides a tangible way to get hold of the C*-envelope.
One starts with a completely isometric representation $\rho$, dilates it to a
maximal representation $\pi$, and $\cenv(\A) = \ca(\pi(\A))$.

In \cite{DKmultdyn}, we provide two views of the C*-envelope of the tensor
algebra $\A(X,\sigma)$.  The first is a rather abstract approach.
Pimsner \cite{Pim} developed a construction of a C*-algebra from a
C*-correspondence, a Hilbert C*-module with a compatible left action,
now known as the Cuntz--Pimsner algebra of the correspondence.
Muhly and Solel \cite{MS2,MS1} developed an extensive theory of an
associated nonself-adjoint \textit{tensor algebra} of a C*-correspondence.
They show that when the left action is faithful, the C*-envelope of the tensor
algebra is the  Cuntz--Pimsner C*-algebra of the C*-correspondence.
Katsura \cite{Ka} extended this theory, defining the Cuntz--Pimsner algebra
for more general left actions which need not be faithful.
The second author and Kribs \cite{KK3}
used Katsura's work to generalize the Muhly--Solel theorem to describe
the C*-envelope of the tensor algebra in full generality.

We show explicitly \cite{DKmultdyn} that the tensor algebra of a
multivariable dynamical system is the tensor algebra of a naturally
associated C*-correspondence.
Consequently, by the results in the previous paragraph, we have a description
of the C*-envelope as a Cuntz--Pimsner algebra.
Unfortunately, because this algebra is a quotient of the
Cuntz--Toeplitz algebra by Katsura's ideal, this is not a very concrete description.
We were looking for something more explicit.

The first step, carried out in \cite{DKmultdyn}, is to describe the maximal
dilations of the orbit representations.
Notice that if $x = \sigma_i(y)$, then the orbit representation $\pi_x$ can be obtained
as the restriction of $\pi_y$ to an invariant subspace.  Hence $\pi_x$ dilates to $\pi_y$.
One can repeat this procedure, and when the system is surjective, construct
an infinite chain of orbit representations, each being a dilation of the previous one.
The inductive limit of this procedure yields a family of maximal representations.
When the system is not surjective, this procedure stops if we arrive at a point $y$
which is not in the range of any map.  It turns out that the orbit representation of
such a point is also maximal.
Since the direct sum of all orbit representations is completely isometric on
the tensor algebra, it follows that the direct sum of all of these maximal dilations
is also completely isometric.  Hence the C*-envelope is given as the algebra
generated by this large representation.
This still is not very explicit, so we seek to develop this some more.

In the case $n=1$ when $\sigma$ is surjective and $X$ is compact,
this was accomplished by Peters \cite{Pet3}.
The idea is to take the projective limit $\tilde X$ of the system $(X,\sigma)$
\[
 X \overset{\sigma}{\longleftarrow}
 X \overset{\sigma}{\longleftarrow}
 X \overset{\sigma}{\longleftarrow}
 \dots \longleftarrow \tilde X .
\]
There is canonical projection $p$ of $Y$ onto $X$, and a map $\tilde\sigma$
of $Y$ onto itself such that $p\tilde\sigma = \sigma p$.  Moreover, $\tilde\sigma$
is always a homeomorphism.
Consequently, one can form the C*-crossed product $\rC(\tilde X) \times_{\tilde\sigma} \bZ$.
There is a natural injection of $\rC(X) \times_\sigma\bZ^+$ into this algebra
by sending $f$ to $f\circ p$ and sending $\fs$ to the canonical unitary of the
crossed product.  Peters shows that this map is a complete isometry, and
that the image generates the crossed product as a C*-algebra.
Then with a bit more work, one obtains:

\begin{thm}[Peters \cite{Pet3}]
Let $X$ be a compact Hausdorff space, and let $\sigma$ be a surjective
continuous map of $X$ onto itself.  Construct $(\tilde X,\tilde\sigma)$ as above.
Then
\[ \cenv(\A(X,\sigma)) \simeq \rC(\tilde X) \times_{\tilde\sigma} \bZ .\]
\end{thm}

The first author and Jean Roydor \cite{DavR} extended Peters' construction to the
multivariable setting.
First assume that $(X,\sigma)$ is surjective in the sense that
$X = \bigcup_{i=1}^n \sigma_i(X)$.
One can again construct a projective limit system.
An infinite tail is an infinite sequence $\bi \in \bn^\bN$, where $\bn = \{1,\dots,n)$;
say $\bi = (i_0,i_1,\dots)$.
One considers the set $\tilde X$ of all pairs $(\bi,\bx) \in \bn^\bN \times X^\bN$
such that $\sigma_{i_s}(x_{s+1}) = x_s$.
There is again a natural map $p(\bi,\bx) = x_0$ of $\tilde X$ onto $X$ and
maps $\tilde\sigma_j(\bi,\bx) = \big( (j,\bi), (\sigma_j(x_0),\bx) \big)$ that
satisfy $p\tilde\sigma_j = \sigma_j p$ for $1 \le j \le n$.
These maps are no longer homeomorphisms.  However the range $\tilde X_j$
of $\tilde\sigma_j$ consists of all points $(\bi,\bx)$ such that $i_0=j$.
These are pairwise disjoint clopen sets, and $\tilde\sigma_j$ is a homeomorphism
of $\tilde X$ onto $\tilde X_j$.
The inverse map $\tau$ given by $\tau|_{\tilde X_j} = \tilde\sigma_j^{-1}$ is a
local homeomorphism.
The tensor algebra $\A(X,\sigma)$ imbeds completely isometrically into
$\A(\tilde X, \tilde\sigma)$, and they have the same C*-envelope.

This leads to a more concrete description of the C*-envelope, because the
new system is much simpler to handle.  One description is that this is the
groupoid C*-algebra of the system $(\tilde X,\tau)$ in the sense of Deaconu \cite{Dea}.
Another is that it is the crossed product of a certain inductive limit $\fB$ of
homogeneous C*-algebras by an endomorphism $\alpha$.

\begin{thm}[Davidson--Roydor \cite{DavR}]
Let $X$ be a locally compact Hausdorff space, and let $\sigma_1,\dots,\sigma_n$
be proper maps of $X$ into itself such that $X = \bigcup_{i=1}^n \sigma_i(X)$.
Construct $(\tilde X,\tilde\sigma)$ as above.
Then
\[ \cenv(\A(X,\sigma)) \simeq \cenv(\A(\tilde X,\tilde \sigma))
\simeq \ca( \tilde X, \tau) \simeq \fB \times_\alpha \bZ^+ .\]
\end{thm}

When $(X,\sigma)$ is not surjective, there is a well-known technique from
graph algebras of ``adding a tail''.
Let $U = X \setminus \bigcup_{i=1}^n \sigma(X)$.
Form $X^T = X \cup T$ where $T = \{(u,k) :  u \in \ol{U},\ k < 0 \}$.
Extend $\sigma_i$ to maps $\sigma^T_i$ by setting
\[ \sigma_i^T(u,k) = (u,k+1) \FOR k < -1, \qand \sigma_i^T(u,-1) = u .\]
It is shown that the natural imbedding of $\A(X,\sigma)$ into $(X^T,\sigma^T)$
is a completely isometric isomorphism.  Moreover, the C*-envelope of
$\A(X,\sigma)$ is a full corner of $\cenv(\A(X^T,\sigma^T))$.

One consequence is a characterization of when the C*-envelope is simple.
When $X$ is compact, the system $(X,\sigma)$ is called \textit{minimal}
if there are no proper closed $\sigma$-invariant subsets of $X$.

\begin{thm}[Davidson--Roydor \cite{DavR}]
Let $(X,\sigma)$ be a compact multivariable dynamical system. Then
$\cenv(\A(X,\sigma))$ is simple if and only if $(X,\sigma)$ is minimal.
\end{thm}

%%%%%%%%%%%%%%%%%%%%%%%%%%%%%%%%

%%%%%%%%%%%%%%%%%%%%%%%%%%%%%%%%
\end{document}